\documentclass[12pt,reqno]{amsart} 

\usepackage{graphicx,amssymb,colordvi,textcomp,latexsym,mathtools,enumitem}
\usepackage[hidelinks]{hyperref}
\usepackage{verbatim}
\usepackage[noend]{algpseudocode}

\def\ddefloop#1{\ifx\ddefloop#1\else\ddef{#1}\expandafter\ddefloop\fi}
\def\ddef#1{\expandafter\def\csname 
	bb#1\endcsname{\ensuremath{\mathbb{#1}}}}
\ddefloop ABCDEFGHIJKLMNOPQRSTUVWXYZ\ddefloop
\def\ddef#1{\expandafter\def\csname 
	#1\endcsname{\ensuremath{\mathbf{#1}}}}
\ddefloop ABCDEFGHIJKLMNOPQRSTUVWXYZ\ddefloop
\def\ddef#1{\expandafter\def\csname 
	#1\endcsname{\ensuremath{\mathbf{#1}}}}
\ddefloop abcdefghijklmnopqrstuvwxyz\ddefloop
\def\ddefloop#1{\ifx\ddefloop#1\else\ddef{#1}\expandafter\ddefloop\fi}
\def\ddef#1{\expandafter\def\csname 
	b#1\endcsname{\ensuremath{\mathbb{#1}}}}
\ddefloop ABCDEFGHIJKLMNOPQRSTUVWXYZ\ddefloop
\def\ddef#1{\expandafter\def\csname 
	c#1\endcsname{\ensuremath{\mathcal{#1}}}}
\ddefloop ABCDEFGHIJKLMNOPQRSTUVWXYZ\ddefloop
\def\ddef#1{\expandafter\def\csname 
	f#1\endcsname{\ensuremath{\mathfrak{#1}}}}
\ddefloop ABCDEFGHIJKLMNOPQRSTUVWXYZabcdefghjklmnopqrstuvwxyz\ddefloop
\def\ddef#1{\expandafter\def\csname 
	h#1\endcsname{\ensuremath{\widehat{#1}}}}
\ddefloop ABCDEFGHIJKLMNOPQRSTUVWXYZ\ddefloop
\def\ddef#1{\expandafter\def\csname 
	hc#1\endcsname{\ensuremath{\widehat{\mathcal{#1}}}}}
\ddefloop ABCDEFGHIJKLMNOPQRSTUVWXYZ\ddefloop

\renewenvironment{proof}{\par\noindent{\bf Proof\ }}{\hfill\BlackBox\\[2mm]}

\newcommand{\tens}{\otimes{}}
\newcommand{\defeq}{\coloneqq}
\DeclarePairedDelimiter{\brk}{[}{]}

\DeclarePairedDelimiter{\prn}{(}{)}
\DeclarePairedDelimiter{\nrm}{\|}{\|}
\DeclarePairedDelimiter{\inner}{\langle}{\rangle}

\usepackage{mathtools}
\usepackage{enumitem}

\newcommand{\symp}{S}	
\newcommand{\sympn}{\symp^n(\reals^d)}
\newcommand{\iso}[1]{(S_k\times S_d)_{#1}}
\newcommand{\reals}{\mathbb{R}}

\newcommand{\Sym}{{\rm Sym}}
\newcommand{\sym}{{\rm sym}}
\newcommand{\ploss}{{\cL}} 
\newcommand{\ibr}[1]{{[#1]}}
\newcommand{\GLG}[1]{{{\text{GL}(#1,\real)}}}

\usepackage{chngcntr}

\theoremstyle{definition}  

\newtheorem{lemma}{Lemma}

\newtheorem{prop}{Proposition}

\theoremstyle{plain}
\newtheorem{rem}{Remark}

\usepackage{prettyref}
\newcommand{\pref}[1]{\prettyref{#1}}

\newcommand{\savehyperref}[2]{\texorpdfstring{\hyperref[#1]{#2}}{#2}}
\newrefformat{eq}{\savehyperref{#1}{\textup{(\ref*{#1})}}}
\newrefformat{eqn}{\savehyperref{#1}{Equation~\ref*{#1}}}
\newrefformat{lem}{\savehyperref{#1}{Lemma~\ref*{#1}}}
\newrefformat{obs}{\savehyperref{#1}{Observation~\ref*{#1}}}
\newrefformat{def}{\savehyperref{#1}{Definition~\ref*{#1}}}
\newrefformat{thm}{\savehyperref{#1}{Theorem~\ref*{#1}}}
\newrefformat{corr}{\savehyperref{#1}{Corollary~\ref*{#1}}}
\newrefformat{sec}{\savehyperref{#1}{Section~\ref*{#1}}}
\newrefformat{app}{\savehyperref{#1}{Appendix~\ref*{#1}}}
\newrefformat{ass}{\savehyperref{#1}{Assumption~\ref*{#1}}}
\newrefformat{ex}{\savehyperref{#1}{Example~\ref*{#1}}}
\newrefformat{fig}{\savehyperref{#1}{Figure~\ref*{#1}}}
\newrefformat{alg}{\savehyperref{#1}{Algorithm~\ref*{#1}}}
\newrefformat{rem}{\savehyperref{#1}{Remark~\ref*{#1}}}
\newrefformat{conj}{\savehyperref{#1}{Conjecture~\ref*{#1}}}
\newrefformat{prop}{\savehyperref{#1}{Proposition~\ref*{#1}}}
\newrefformat{proto}{\savehyperref{#1}{Protocol~\ref*{#1}}}
\newrefformat{prob}{\savehyperref{#1}{Problem~\ref*{#1}}}
\newrefformat{claim}{\savehyperref{#1}{Claim~\ref*{#1}}}
\newrefformat{fact}{\savehyperref{#1}{Fact~\ref*{#1}}}
\newrefformat{exam}{\savehyperref{#1}{Example~\ref*{#1}}}
\newrefformat{exams}{\savehyperref{#1}{Examples~\ref*{#1}}}
\newrefformat{prob}{\savehyperref{#1}{(\ref*{#1})}}

\newtheorem{exam}[lemma]{Example}
\newtheorem{exams}[lemma]{Example}

\newcommand{\scalepar}{0.4}

\usepackage{setspace}

\newcommand{\RR}{\mathbb{R}}
\newcommand{\EE}{\mathbb{E}}

\newcommand{\PP}{\mathbb{P}}

\newcommand{\NN}{\mathbb{N}}

\renewcommand{\ker}{{\kappa}}

\newcommand{\pint}{\mbox{$\mathbb{N}$}}

\usepackage{chngcntr}
\counterwithout{equation}{subsection}
\counterwithout{equation}{section} 
\renewenvironment{proof}{\par\noindent{\bf Proof\ }}{\hfill\\[2mm]}

\begin{document}

\author{Yossi Arjevani, Joan Bruna, Michael Field, Joe Kileel, Matthew Trager \and Francis Williams}
\address{Yossi Arjevani, The Hebrew University}
\address{Joan Bruna, New York University}
\address{Michael Field,  UC Santa Barbara}
\address{Joe Kileel, University of Texas at Austin}
\address{Matthew Trager, Amazon (part of the work was performed prior to joining Amazon)}
\address{Francis Williams, New York University}

\title[Symmetry Breaking in Symmetric Tensor Decomposition]{Symmetry Breaking in Symmetric\\Tensor Decomposition}

\date{}

\begin{abstract}
In this note, we consider the highly nonconvex optimization problem associated with computing the rank decomposition of symmetric tensors. We formulate the invariance properties of the loss function and show that critical points detected by standard gradient based methods are \emph{symmetry breaking} with respect to the target tensor. The phenomena, seen for different choices of target tensors and norms, make possible the use of recently developed analytic and algebraic tools for studying nonconvex optimization landscapes exhibiting symmetry breaking phenomena of similar nature.
\end{abstract}
\maketitle
\markleft{}

\section{Introduction}

We consider the problem of approximating a symmetric tensor as a sum of 
rank-1 symmetric tensors. Concretely, given an order $n$ symmetric tensor $T\in 
{\rm Sym}^n(\RR^d)$ and rank $k\in\NN$, we study the nonconvex optimization problem
\begin{align} \label{prob:opt}
	\min_{W\in M(k,d)}\ploss_\star(W)  \defeq  
	\nrm*{ \sum_{i=1}^k \w_i^{\tens 
			n} - T}_\star^2,
\end{align}
with $M(k,d)$ denoting the space of $k\times d$ matrices, $\nrm{\cdot}_\star$ a tensor norm, and $\w_i$ the $i$-th row of $W$. We emphasize odd $n$ where (\ref{prob:opt}) is equivalent to standard tensor rank decomposition, also known as the real symmetric \emph{canonical polyadic decomposition (CPD)}. However, methods are quite general.

The problem of finding a tensor approximation of bounded rank arises naturally in various scientific fields, including machine learning, biomedical engineering and psychometrics, see \cite{comon2002tensor,comon2006blind,de2000multilinear, kolda2009tensor, sidiropoulos2000parallel,smilde2005multi, landsberg2011tensors} and references therein for applications. The associated optimization problem is highly nonconvex and exhibits a variety of saddles and spurious (i.e., non-global local) minima that can cause a complete failure of gradient-based optimization methods. It is therefore of interest to study and characterize geometric obstructions of the this nature that exist in the associated loss landscape. 

\begin{figure}
	\centering \includegraphics[scale=\scalepar]{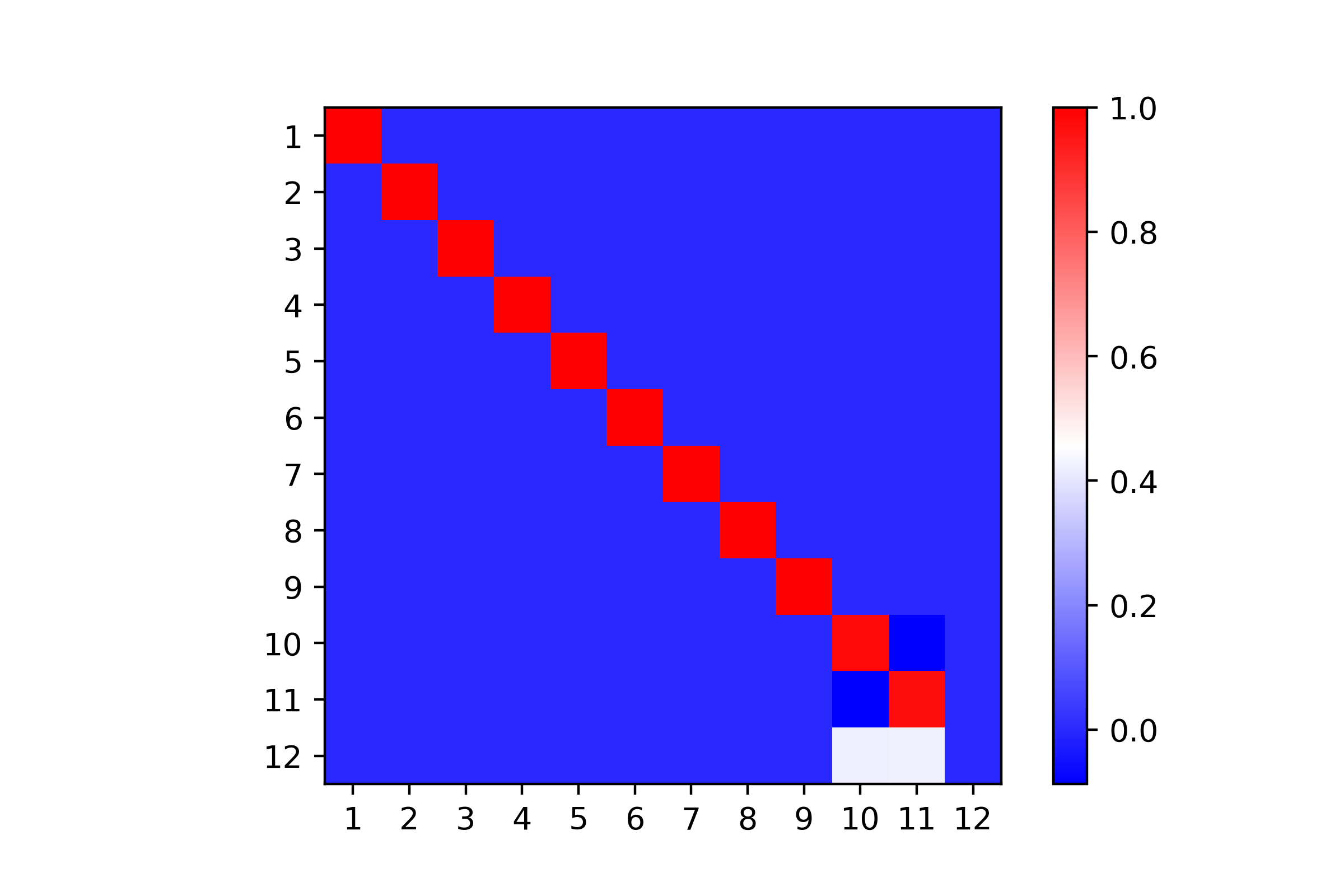}
	\includegraphics[scale=\scalepar]{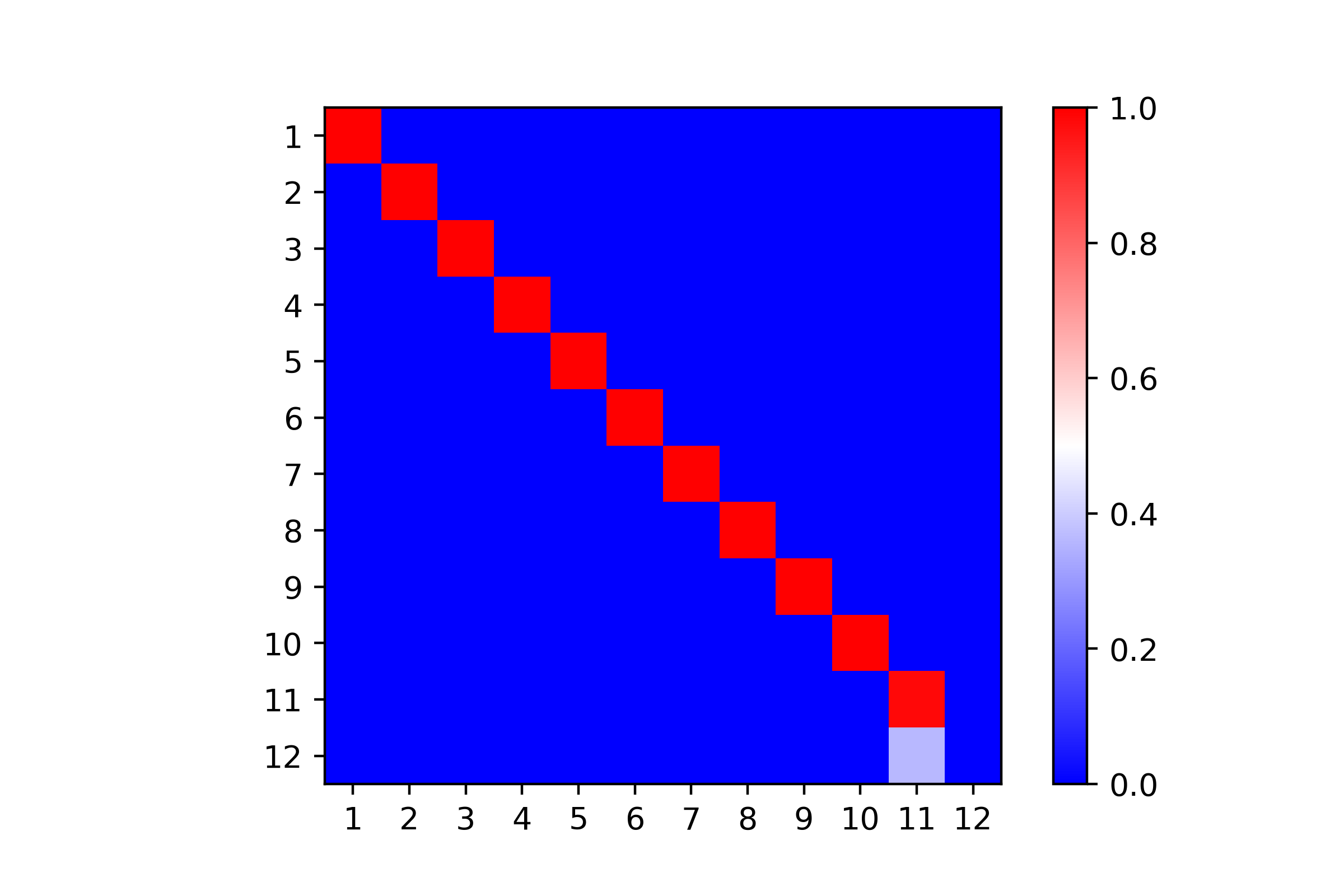}
	\caption{Spurious minima of optimization problem \pref{prob:opt} under the Frobenius norm and problem parameters $d=k=12$, order $3$ symmetric tensors
and target tensor $T= 	\sum_{i=1}^k \e_i^{\tens 3}$. (Left) a $12\times 
		12$-spurious minimum fixed by the group of row and column permutations $\Delta (S_{9}\times S_2 \times S_1)\subseteq S_{12}\times S_{12}$. (Right) a $12\times 12$-spurious minimum fixed by $\Delta 
		(S_{10}\times S_1\times S_1)$. }
	\label{fig:r3_a}
\end{figure}

In this note, we show that, empirically, critical points of $\ploss_\star$ are \emph{symmetry breaking}, a property defined in terms of \emph{isotropy groups}, i.e., the group of all row and column permutations, denoted by $(S_k\times S_d)_{C}\subseteq S_k\times S_d$, that fix the weights of a given critical point $C \in M(k,d)$. We found that critical points detected by gradient-based methods  $\iso{C}$ are typically large and conform with the structure of the target tensor~$T$ (and so lie on low-dimensional subspaces), see \pref{fig:r3_a}. The presence of symmetry breaking phenomena makes possible the derivation of sharp analytic estimates for families of critical points, their loss and their Hessian spectrum \cite{arjevani2023symmetry}.

Our contributions in order of appearance may be stated as follows.  In \pref{sec:symmetry}, invariance properties of $\ploss_\star$ are studied for different choices of inner products and target tensors, emphasizing the  dependence of the invariance properties on the target tensor. In \pref{sec:numerical_results}, numerical results are provided, indicating various isotropy groups of critical points that are found numerically. Results are given for the Frobenius inner product of order 3 and order 5 tensors, and for the inner product associated with the standard Gaussian distribution. Relevant background from multilinear algebra and group action is briefly reviewed in \pref{sec:preliminaries}.\\

Next, we relate our results to the existing literature.

\subsubsection{Symmetric tensor decomposition.}
We briefly mention current algorithms for symmetric tensor 
decomposition.  A straightforward yet practical method is based on direct 
first-order optimization of $\mathcal{L}$ in \eqref{prob:opt}; see 
\cite{kolda2015numerical} and the Matlab implementations \cite{tensorlab, 
	tensortoolbox}.  A more computationally intensive although provable method 
(assuming $T$ is rank $r$ and $r = O(n^{\frac{d-1}{2}})$) was provided in 
\cite{nie2017generating}, based on an algebraic construction called generating 
polynomials. For $n=3$ and $r \leq d$, a classic but theoretically convenient 
method was derived from simultaneous diagonalization of matrix slices 
\cite{harshman70} .   For $n=4$ and $r = O(n^2)$, the work \cite{de2007fourth} 
presented a provable algorithm using matrix eigendecompositions, which was 
robustified using ideas from the sums-of-squares hierarchy in 
\cite{hopkins2019robust}.  In \cite{kileel2019subspace}, a tensor power method 
was used, constructed from a matrix flattening of $T$, to find the components 
$\w_i$ sequentially.  Also, \cite{sherman2020estimating} showed how to 
implement direct optimization of \eqref{prob:opt} in an efficient manner for 
moment tensors $T$ in an online setting.

\subsubsection{Symmetry breaking in nonconvex loss landscapes.}
It has been recently found that spurious minima occurring for a version of optimization problem (\ref{prob:opt}) defined by use of the ReLU activation {break the symmetry} of global minima \cite{ArjevaniField2019} (a formal exposition shall be given later in \pref{sec:symmetry}). Here, we report phenomena of symmetry breaking occurring for symmetric tensor decomposition problems. These were later used in \cite{ArjevaniField2022} 
to study families of critical points building on methods developed in a line of work concerning phenomena of symmetry breaking, which we now survey. In 
\cite{arjevani2021symmetry}, path-based techniques are introduced, 
allowing the construction of infinite families of critical points 
for ReLU two-layer networks using Puiseux series. In 
\cite{arjevani2020analytic}, results from the representation 
theory of the symmetric group are used to obtain precise analytic 
estimates on the Hessian spectrum. In~\cite{arjevani2021analytic}, it is 
shown that certain families of saddles transform into spurious minima 
at a fractional dimensionality. In addition, Hessian spectra at 
spurious minima are shown to coincide with that of global minima 
modulo $O(d^{-{1}/{2}})$-terms. In 
\cite{arjevani2022annihilation}, it is proved that adding 
neurons can turn symmetric spurious minima into saddles. In 
\cite{arjevani2022equivariant}, generic $S_d$-equivariant 
steady-state bifurcation is studied, emphasizing irreducible representations along which 
spurious minima may be created and annihilated. In \cite{arjevani2023hidden}, it is shown that the way subspaces \emph{invariant} to the action of subgroups of $S_d$ are arranged relative to ones \emph{fixed} by the action determines the admissible types of structure and symmetry of curves along which $\ploss$ is minimized and maximized.

\section{Preliminaries}\label{sec:preliminaries}
We  review relevant background material from multilinear algebra and group actions needed for a formal study of symmetry breaking.

\subsection{Tensor preliminaries} 

A tensor of order $n$ is an element of the tensor 
product of $n$ vector spaces $V_1 \otimes \ldots \otimes V_n$. Assuming $V_i$ are real vector spaces, we may choose a base for each factor $V_i$ and so identify a tensor with a multi-dimensional array in $\RR^{d_1 	\times \ldots \times d_n}$ with $d_i = \dim(V_i)$. We write $T_{i_1,\ldots,i_n}$ for the $(i_1,\ldots,i_n)$-th 
coordinate of $T$. Given vectors $\v_i \in \RR^{d_i}$, $i=1,\ldots,n$, 
we write $\v_1 \otimes \ldots \otimes \v_n$ for the outer product of these 
vectors, that is, the element 
in $\RR^{d_1} \otimes \ldots \otimes \RR^{d_n}$ such that  $(\v_1 \otimes 
\ldots \otimes \v_n)_{i_1,\ldots,i_n} = ({v_1})_ {i_1} \ldots (v_n)_{i_n}$. 
Also, we write $\v^{\otimes n} := \v \otimes \ldots \otimes \v$ ($n$-times). The 
Frobenius inner product of two tensors $T, S$ of the same shape is defined to be
$\langle T, S \rangle_F = \sum_{i_1,\ldots,i_n} T_{i_1\ldots i_n} S_{i_1 
	\ldots i_n}$. We will also consider other inner products of tensors. For the Frobenius inner product,
\begin{equation}\label{eq:frobenius_vector_product}
	\langle \v_1 \otimes
	\ldots \otimes \v_n, \w_1 \otimes \ldots \otimes \w_n \rangle_F = \prod_{i=1}^n \langle \v_i, \w_i  \rangle
\end{equation}
for all $\v_i, \w_i \in \mathbb{R}^{d_i}, i=1, \ldots, n$.  In particular, 
$\langle \v^{\otimes n}, \w^{\otimes n}\rangle_F = \langle \v,  
\w\rangle^n$.

A tensor $T \in (\RR^{d})^{\tens n}$ is \emph{symmetric} if it 
is invariant under permutation of indices, that is, if 
$T_{i_1,\ldots,i_n} = 
T_{i_{\sigma(1)},\ldots,i_{\sigma(n)}}$ for any
permutation $\sigma \in S_n$. The space of symmetric 
tensors of order $n$ on $\RR^d$, denoted by ${\rm \symp}^n(\RR^d)$, is 
$\binom{d+n-1}{n}$-dimensional and is isomorphic to the space of 
homogeneous polynomials of degree $n$ in $d$ variables. An natural 
isomorphism is given by the map 
\begin{align}\label{sp_id}
	S\mapsto P_S,~~P_S(x) \defeq \inner{ S, 	x^ {\otimes r}}_F.
\end{align}
When $n=2$, this is the usual correspondence between 
symmetric matrices and quadratic forms. We consider two 
choices of inner products for ${\rm \symp}^n(\RR^d)$: the restriction of the Frobenius inner product 
$\inner{\cdot,\cdot}_F$ to $\symp^n(\reals^d)$, and 
\begin{align}\label{data_ip}
	\inner{S, T}_{\cD} \defeq \mathbb{E}_{\x \sim \mathcal{D}} 
	[\inner{S, 
		\x^{\otimes 
			n}}_F \inner {T, 	\x^{\otimes n}}_F],
\end{align}	
where $\cD$ is a distribution on $\reals^d$ dominating the Lebesgue 
measure and chosen so that the above expectation is finite for any $S,T\in 
\symp^n(\reals^d)$.

\begin{prop}[Notation and assumptions as above]\label{prop:data_ip} 
	The bivariate function $\inner{\cdot, 
		\cdot}_{\cD}:\symp^n(\reals^d)\times\symp^n(\reals^d)\to\reals$ defines an 
	inner product on $\sympn$.
\end{prop}
Note that over $\reals^{\tens d}$, $\inner{\cdot,\cdot}_\cD$ is \emph{not} positive definite, and so the  explicit restriction to $\sympn$ is necessary. 
\proof 
It is easy to verify that $\inner{\cdot,\cdot}_\cD$ 
is a symmetric bilinear form. To prove that $\inner{\cdot,\cdot}_\cD$ is positive 
definite, observe that for any $S\in\symp^n(\reals^d)$,
\begin{align}\label{ten_psd}
	\inner{S, S}_{\mathcal D} &= \mathbb E_{\x \sim 
		\mathcal D}\brk{ (\inner{S, \x^{\otimes n}}_F)^2} \ge 0.
\end{align}
Now, let $S$ be such that 
(\ref{ten_psd})~holds with equality, and assume by way of contradiction that there 
exists $\x_0\in\reals^d$ such that $P_S(\x_0)\defeq 
\inner{S, \x_0^{\otimes n}}_F \neq 
0$. By continuity, there exists an open set $U\in\reals^d$ such that 
$P^2_S(\x)>P_S^2(\x_0)/2 >0$ on $U$. Therefore, by the law of total 
expectation,
\begin{align}
	\inner{S, S}_{\mathcal D} \ge \EE[P^2_S(\x)~|~U]\PP[U]\ge \PP[U]P_S^2(\x_0)/2>0,
\end{align}
contradicting equality in (\ref{ten_psd}). Thus, $\inner{S, 
	\x^{\otimes n}}_F$ is a homogeneous polynomial vanishing everywhere on 
$\reals^d$, and so (necessarily over a field of characteristic zero) $S=0$, 
concluding the proof. \qed\\
\\
We shall be particularly interested in the \emph{cubic-Gaussian} inner product defined by setting $\cD=\cN(0, I_d)$ and $n=3$. We denote the corresponding loss function (over $\sym^3(\reals^d)$) by  $\ploss_{\cN}$.

In terms of the coefficients of $P_S(x) = \sum 
a_{p_1,\ldots,p_d} x_1^{p_1}\ldots x_d^{p_d}$ and $P_T(x) = \sum 
b_{p_1,\ldots,p_d} x_1^{p_1}\ldots x_d^{p_d}$, the Frobenius inner product 
on $\symp^n(\RR^d)$ is given by 
\begin{equation}\label{eq:frobenius_bombieri}
	\inner{S, T}_F = \sum_{p_1+ \ldots + p_d = n} \binom{n}{p_1, \ldots,
		p_d} a_{p_1,\ldots,p_d}	b_{p_1,\ldots,p_d}.
\end{equation}
The inner product $\inner{\cdot,\cdot}_\cD$ reads
\begin{equation}\label{eq:distributional_product}
	\inner{T, S}_{\mathcal D} = \sum_{\substack{p_1+ \ldots + p_d = n\\
			q_1+ \ldots + q_d = n}} \mathbb E[x_1^
	{p_1 + q_1}\ldots x_d^{p_d + q_d}] a_{p_1,\ldots,p_d} 
	b_{q_1,\ldots,q_d}.
\end{equation}
For the standard multivariate normal distribution $\cD= 
\cN(0,I_d)$, the above may be given in an explicit form,
\begin{align}
	&\inner{S, T}_{\cN} = \sum h_{\p+\q} a_
	{p_1,\ldots,p_d} b_{p_1,\ldots,p_d},
	\\
	&\mbox{ where } h_{\r} := 
	\begin{cases} 0 & \mbox{if } r_i \mbox{ is odd for some } i \in
		\{1,\ldots,d\},\\
		\prod_{i=1}^d {(r_i-1) !!} & \mbox{otherwise. }
	\end{cases}\nonumber
\end{align}

Comparing \eqref{eq:frobenius_bombieri} and 
	\eqref{eq:distributional_product}, it is seen that no data distribution induces 
	the Frobenius product. 
	Indeed, in~\eqref{eq:distributional_product}, whether a term 
	$a_{p_1,\ldots,p_d}$ is multiplied by $b_{q_1,\ldots,q_d}$ depends only on 
	${p_1 + q_1},\ldots,{p_d
		+ q_d}$, and so in particular if the coefficient of $a_ 
	{p_1,\ldots,p_d}
	b_{q_1,\ldots,q_d}$ is non-zero then, assuming, e.g., $p_1 > 0$, so is the 
	coefficient corresponding to $a_{p_1-1,\ldots,p_d} b_{q_1+1,\ldots,q_d}$. 
	In the Frobenius product however the coefficient of $a_ {p_1,\ldots,p_d} b_ 
	{q_1,\ldots,q_d}$ is non-zero if and only if $(p_1,\ldots,p_d) = 
	(q_1,\ldots,q_d)$. 

\begin{rem}
In the language of homogeneous polynomials, by considering problem 
\eqref{prob:opt} with structured targets $T$, we are studying the 
\textit{Waring decompositions} of partially symmetric homogeneous polynomials  
\cite{landsberg2011tensors}. 
\end{rem}

The data-dependent inner product defined in (\ref{data_ip}) can also be expressed in terms a similarity measure between vectors in $\reals^d$, often referred to as a \emph{kernel}. Given $\w,\v\in\reals^d$, we define 
\begin{align} \label{ker_cd}
	\ker_\cD(\w,\v) \defeq \bE_{\x \sim \cD} [\rho(\inner{\w,\x})\rho(\inner{\v,\x})],
\end{align}
with $\cD$ denoting a distribution over $\reals^d$ and $\rho:\reals\to\reals$ a measurable function.
For example, the inner product defined by the choice of 
$\rho(z)=\max\{z,0\}$, the ReLU activation function, and $\cD=\cN(0,I_d)$ is 
used in the study of two-layer ReLU neural networks 
\cite{ArjevaniField2019}.

\begin{prop}[Notation and assumptions as above]
	\label{prop:tensor_ip}
Let $\rho(z)=z^n$. If $S,T\in \symp^n(\reals^d)$ are given by 
$S= \sum_{i=1}^k \alpha_i\w_i^{\tens n}$ and $T = \sum_{j=1}^h\beta_j\v_j^{\tens n}$, with $\w_i$ (resp. $\v_i$) denoting $k$ (resp. $h$) $d$-dimensional vectors, then 
\begin{align}\label{ip_in_ker}
\inner{S, T}_{\cD} = \sum_{i=1}^{k} \sum_{j=1}^{h} \alpha_i\beta_j\ker_\cD(\w_i, \v_j).
\end{align}		 
\end{prop}
The expression of $\inner{S,T}_\cD$ in terms of $\ker_\cD$ comes handy when studying 
the invariance properties of the loss function, see \pref{sec:symmetry}.
\proof
The result follows by a direct computation.
\qed\\
\\
We refer to the data-dependent kernel corresponding to $\cD  = \cN(0,I_d)$ and $\rho(z)=z^3$ as the \emph{cubic-Gaussian} kernel, given in explicit terms by
\begin{align}\label{kernel:cubic_gauss}
	\ker_{\cN,3}(\w,\v) &= 
	{6}\inner{\w,\v}^3 + 9\|\w\|_2^2\|\v\|_2^2\inner{\v,\w}.
\end{align}
This formula is a direct consequence of Isserlis's Theorem~\cite{isserlis1918formula}, stating that for any $2r$ $d$-dimensional vectors $\v_1,\ldots,\v_{2r}\in \RR^{d}$,
\begin{align}
\mathbb E_{\x \sim \mathcal N(0,I_d)} [(\v_1\cdot \x)(\v_2 \cdot \x)&
\dots ( \v_{2r}\cdot \x)] \\&=  \frac{1}{2^r r!}\sum_{\sigma \in S_{2r}} 
(\v_ {\sigma(1)}\cdot \v_ {\sigma(2)}) \ldots (\v_{\sigma({2r-1})}\cdot
\v_{\sigma({2r})}).\nonumber
\end{align}

Finally, we recall the definition of \emph{rank} for symmetric tensors. A symmetric tensor $T \in \Sym^n(\RR^d)$ is said to have  rank-$1$ if $T = \lambda \v^{\otimes n}$ 
for some $\lambda \in \RR \setminus \{0\}$ and $\v \in \RR^d \setminus \{0\}$. More generally, a tensor $T$ has (real symmetric) rank-$r$ if it can be written as a linear combination of $r$ rank-$1$ tensors 
$T = \lambda_1 \v_1^{\otimes n} + \ldots + \lambda_r \v_r^{\otimes n}$, but not 
as a combination of $r-1$ rank-$1$ tensors. For $n=2$, this definition agrees 
with the usual notion of rank for symmetric matrices. We shall only be interested in symmetric tensors of odd order, and so $\lambda_i$ may be absorbed into the vectors $\v_i$ and be dropped altogether.

\subsection{Groups, actions and symmetry}
We start with an example that is used later. Elementary concepts from group theory 
are assumed known. 

\begin{exam} \label{exams:groups}
	Let $\text{GL}(d)$ denote the set of invertible 
	linear maps on $\reals^d$. Under composition, $\text{GL}(d)$ is a 
	group, called the \emph{general linear group}. The \emph{orthogonal group} $\text{O}(d)$ is the 
	subgroup 
	of $\text{GL}(d)$ that preserves Euclidean distances, i.e.,
	$
	\text{O}(d) =   \{A \in \text{GL}(d) : \|Ax\|_2 = 
	\|x\|_2\;\text{for all } x \in \reals^d\}.
	$
	Upon choosing a basis for $\mathbb{R}^d$, both  $\GLG{d}$ and $\text{O}(d)$ can be viewed as groups of 
	invertible $d \times d$ matrices.
\end{exam}

Groups often arise as
\emph{transformations}  
of a set or space, so we are led to the notion of a 
\emph{$G$-space} $X$ where we have 
an \emph{action} of a group $G$
on a set $X$. Formally, a group action is a group homomorphism from $G$ 
to the
group of bijections of $X$.
For example, $S_d$ naturally acts on $[d]:=\{1,\ldots,d\}$ as the group of permutations, while 
both $\textup{GL}(d)$ and $\text{O}(d)$ act on
$\reals^d$ as groups of linear transformations (or matrix 
multiplication). 
\begin{exams} \label{ex:product}
	Our study of the invariance 
	properties of $\ploss_\star$ will rely on the action of the
	product group $S_k \times S_d\subseteq S_{k\times d},~k,d\in 
	\pint$ on the product set
	$\ibr{k} \times \ibr{d}$ defined by
	\begin{align}\label{eq: Gamma-action}
		(\pi,\rho)(i,j) := (\pi^{-1}(i),\rho^{-1}(j)) \,\, \textup{ for all } (\pi, \rho) \in 
		S_k \times S_d, (i,j) \in \ibr{k} \times \ibr{d}.
	\end{align}
	By identifying $(i,j) \in \ibr{k} \times \ibr{d}$ with the entry $(i,j)$-entry in a matrix, this  induces a linear representation on the space $M(k,d)$ of real $k \times 
	d$ matrices $A = \left(A_{ij}\right)$ via
	\begin{equation}
		(\pi,\rho)\left(A_{ij}\right) := \left(A_{\pi^{-1}(i),\rho^{-1}(j)}\right).
	\end{equation}
	Here, $\pi$ acts by permuting the rows of $A$, and $\rho$ by permuting the 
	columns of $A$. In terms of permutation matrices $P_{\pi} \in O(k)$ where 
	\begin{align*}
		(P_\pi)_{ij} = \begin{cases}
			1 & i =\pi(j),\\
			0 & \text{otherwise},
		\end{cases}
	\end{align*}
	and $P_{\rho} \in O(d)$ similarly defined, the action is 
	$(\pi,\rho)(A) = P_{\pi} A P_{\rho}^{\top}$.

\end{exams}

Given a weight matrix $W \in M(k,d)$, the largest 
subgroup of $S_k \times S_d$ fixing $W$ is called the \emph{isotropy} subgroup 
(or \textit{stabilizer})
of $W$ and is used as a means of measuring the symmetry of $W$. For example, the isotropy 
subgroup of $I_d\in M(d,d)$ is the diagonal subgroup $\Delta S_{d}$, where 
$\Delta$ maps a given subgroup $G\subseteq S_d$ to its \textit{diagonal} 
counterpart $\{(g,g) : g \in G\} \subseteq S_d \times S_d$. The \textit{fixed 
	point space} corresponding to a given $G\subseteq S_k \times S_d$ is defined by
\begin{align*}
	M(k,d)^G \defeq \{W\in M(k,d) ~|~ gW = W \,\, \forall h\in G \},
\end{align*} 
forming a linear subspace of $M(k,d)$. Critical points 
with isotropy, e.g., $\Delta S_d$, $\Delta (S_{d-1}\times S_1)$ and $\Delta 
(S_{d-3}\times S_2 \times S_1)$  lie therefore on linear subspaces of fixed dimensionality (respectively, 2, 5 and 11), see  \pref{sec:identity}.

\section{Invariance properties} \label{sec:symmetry}
A real-valued function $f$ with domain $X$ is 
\textit{$G$-invariant} if $f(g\x)=f(\x)$ for all $\x\in X,~g\in G$. Regard $\ploss_\star$ as a function on the $S_k \times S_d$-space $M(k, d)$ (see 
Example~\ref{ex:product}). Permuting the 
order of the summation in \pref{prob:opt} does not change the 
function value, and so $\ploss_\star$ is left $S_k$-invariant, i.e., invariant under 
row permutations of $W$, for any choices of $T$. Additional invariance properties of $\ploss_\star$ are given  by the structure of the target tensor $T$. 

We extend the definition of $\ploss_\star$ so as to make the dependence on the target tensor explicit, writing, by a slight abuse of notation, 
\begin{align} \label{prob:opt_ext}
	\ploss_\star(W;V) \defeq \nrm*{\sum_{i=1}^k \w_i^{\tens n}- \sum_{i=1}^h \v_i^{\tens n}}_\star^2,
\end{align}
with $\w_i$ (resp. $\v_i$) denoting the rows of $W\in M(k,d)$ (resp. $V\in M(h,d)$). Clearly, $\ploss_\star$ is left $S_h$-invariant with respect to the second argument $V$, i.e., invariant to row permutations of $V$.  Under the Frobenius inner product, 
\begin{align} \label{prob:kernel_frob}
	\ploss_F&(W;V) =  \inner{\sum_{i=1}^k \w_i^{\tens n}- \sum_{i=1}^h \v_i^{\tens n}-\sum_{i=1}^k \w_i^{\tens n}- \sum_{i=1}^h \v_i^{\tens n}}_F\nonumber\\
	&=  \sum_{i,j=1}^{k} \inner{\w_i^{\tens n}, \w_j^{\tens n}}_F
	-2\sum_{i=1}^{k} \sum_{j=1}^{h} \inner{\w_i^{\tens n}, \v_j^{\tens n}}_F
	+ \sum_{i,j=1}^{h} \inner{\v_i^{\tens n}, \v_j^{\tens n}}_F\nonumber\\
	&=  \sum_{i,j=1}^{k} \inner{\w_i, \w_j}^n
	-2\sum_{i=1}^{k} \sum_{j=1}^{h} \inner{\w_i, \v_j}^n
	+ \sum_{i,j=1}^{h} \inner{\v_i, \v_j}^n.
\end{align}
The inner product $\inner{\cdot,\cdot}^n$ is invariant to the action of $O(d)$ on $\reals^d$, hence
$\ploss_F(W;V)=\ploss_F(WU;VU)$ for all $U\in O(d)$. By \pref{prop:tensor_ip}, 
\begin{align}\label{loss_in_ker}
\ploss_{\cN}&(W;V) =\nrm*{\sum_{i=1}^k\w_i^{\tens n} - \sum_{j=1}^h 
\v_j^{\tens 	n}}^2_\cN\\
&=
\sum_{i,j=1}^{k} \ker_{\cN,3}(\w_i, \w_j)
-2\sum_{i=1}^{k}  \sum_{j=1}^{h} \ker_{\cN,3}(\w_i, \v_j)
+ \sum_{i,j=1}^{h} \ker_{\cN,3}(\v_i, \v_j). \nonumber
\end{align}
The kernel $\ker_{\cN,3}$ is also invariant to the action of $O(d)$ on $\reals^d$ as seen by (\ref{kernel:cubic_gauss}), and so by the same argument above, $\ploss_\cN(W;V)=\ploss_\cN(WU;VU)$ for all $U\in O(d)$. Thus, in the following, arguments used to derive the invariance properties of $\ploss_\star$ apply equally to $\star = F$ and $\star = \cN$.

Assuming momentarily that $h=d$ and $V = I_d$, we have for all 
$(\pi,\rho)\in S_k\times S_d$:
\begin{align}\label{eq:invariance}
	\ploss_\star(W) &= \ploss_\star(W ; I_d) = \ploss_\star(W ; P_\rho P_\rho^\top )
	= \ploss_\star(W P_\rho ; P_\rho )\nonumber
	\\&= \ploss_\star(P_\pi W P_\rho ; I_d ) = \ploss_\star(P_\pi W P_\rho).
\end{align}
Thus, $\ploss_\star$ is $S_k \times S_d$-invariant when $V=I_d$. More generally, if $V \in M(h,d)$ and $\sigma \in S_h$, $\rho \in S_d$ are such 
that $V P_\rho = P_{\sigma} V$, then
\begin{align}\label{eq:invariance_right}
	\ploss_\star(W) &= \ploss_\star(W ; V) = \ploss_\star(W P_\rho ; V P_\rho) = \ploss_\star(W P_\rho ; 
	P_\sigma V) \nonumber\\&= \ploss_\star(W P_\rho ; V) = \ploss_\star(W P_\rho). 
\end{align}
It is seen that symmetries of the target tensor are reflected by the invariance properties of the loss function. For example, if $V \in M(d,d)$ is a \emph{circulant matrix}, then $V P_\rho = 
P_\rho V$ holds for any \emph{cyclic} permutation $\rho \in S_d$. Thus, 
\eqref{eq:invariance_right} implies that 
$\ploss_\star(W) = \ploss_\star(P_\pi W P_\rho)$ for any permutation $\pi\in S_k$ and any
cyclic permutation $\rho\in S_d$. Symmetry breaking for circulant matrices is discussed in  the next section.

\section{Numerical results}
\label{sec:numerical_results}
\begin{figure}
	\centering \includegraphics[scale=\scalepar]{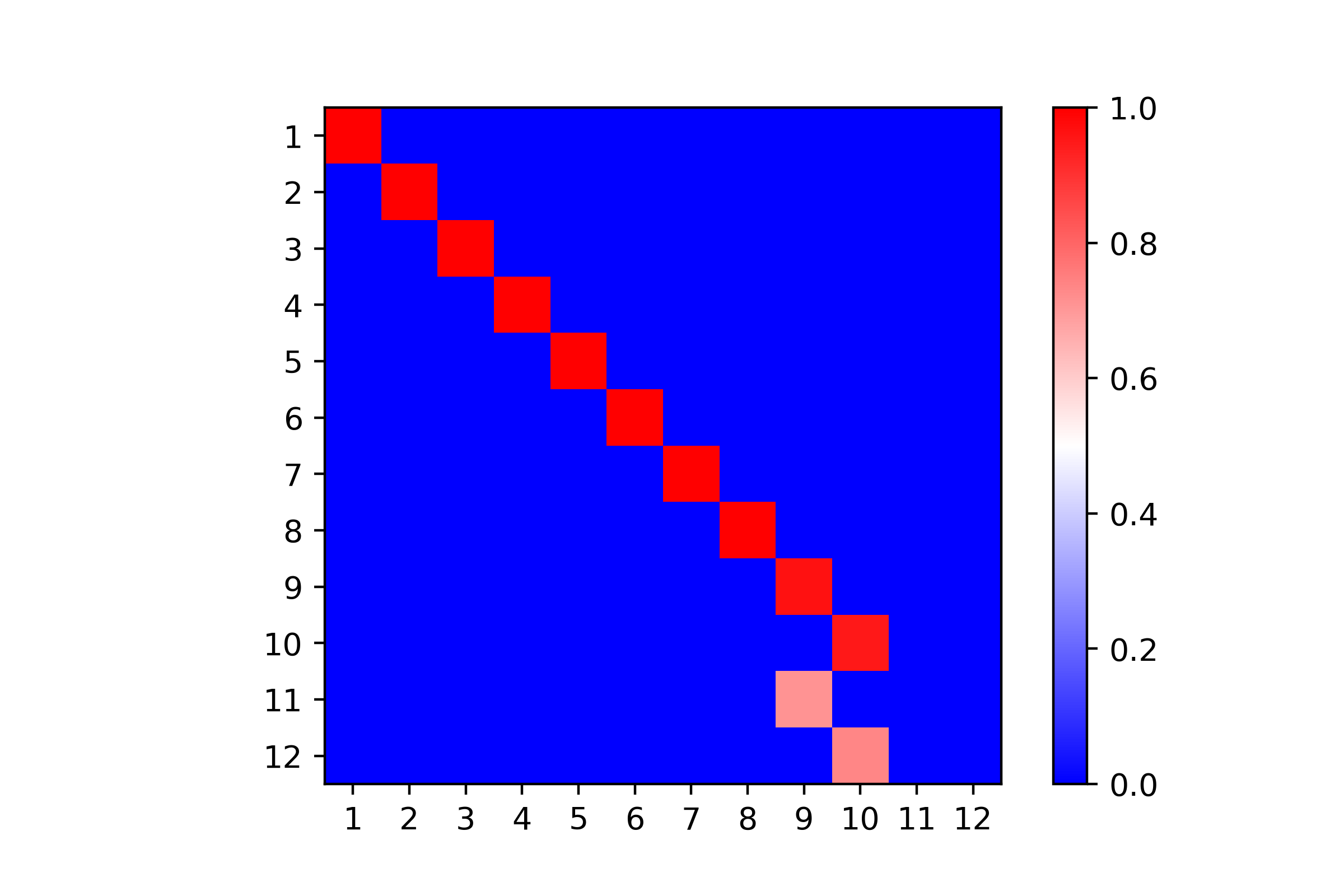}
	\includegraphics[scale=\scalepar]{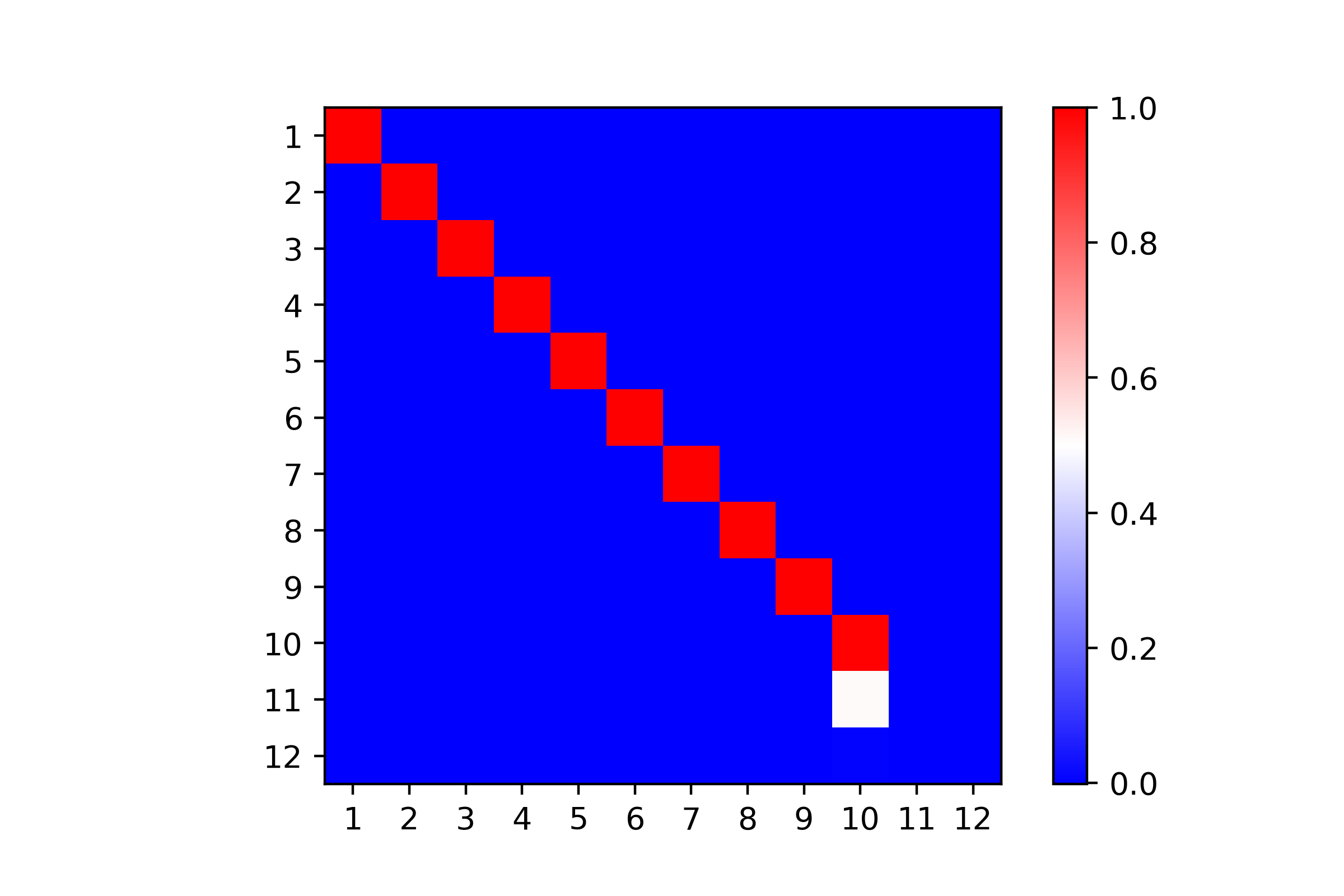}
	\caption{Spurious minima of optimization problem \pref{prob:opt} under the Frobenius norm and
		with problem parameters $d=k=12,~n=5$ and $T= 
		\sum_{i=1}^k \e_i^{\tens 5}$. (Left) a spurious minimum of isotropy 
		$\Delta (S_8\times \inner{(9~10)(11~12)})$. 
		(Right) a spurious 	minimum of isotropy $\Delta (S_9\times S_1^3)$.}
	\label{fig:r5_a}
\end{figure}

\begin{figure}
	\centering \includegraphics[scale=\scalepar]{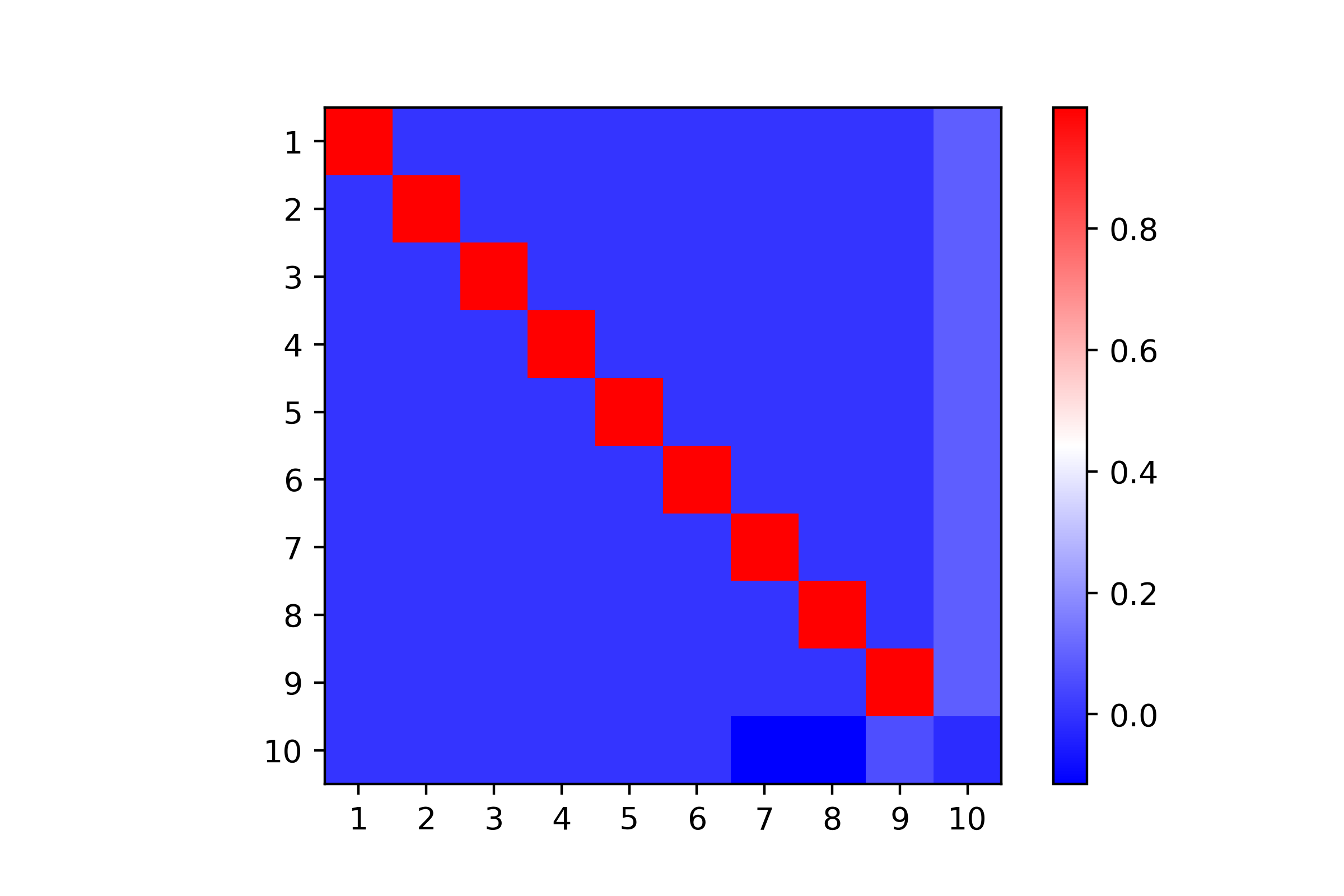}
	\includegraphics[scale=\scalepar]{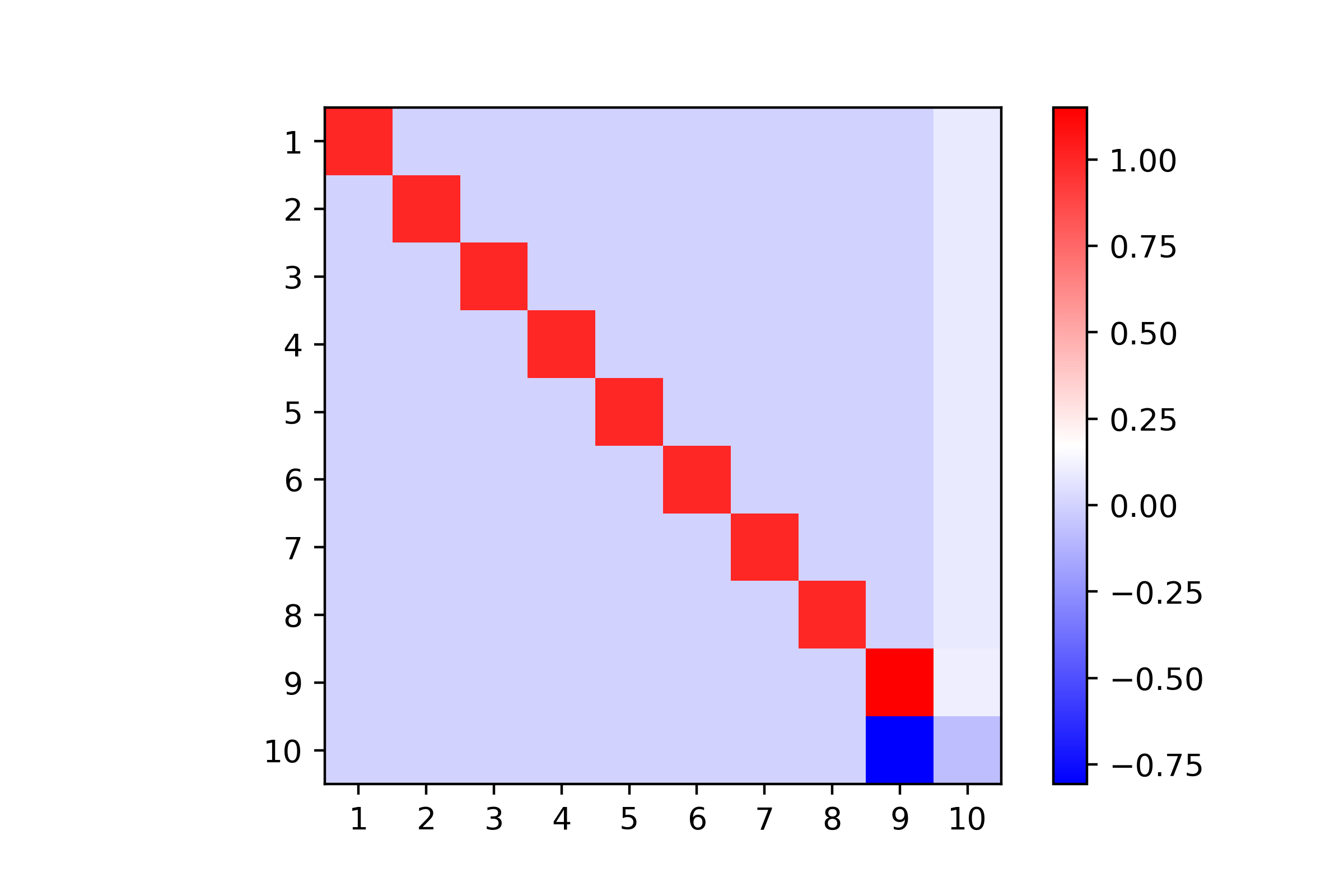}
	\caption{
		Spurious minima of optimization problem \pref{prob:opt} under the cubic-Gaussian  norm and problem parameters $d=k=10,~n=3$ and $T= 	\sum_{i=1}^d \e_i^{\tens 3}$. (Left) a spurious minimum of 
		isotropy $\Delta (S_{6}\times S_2 \times  S_1^2)$.  (Right) 
		a spurious minimum of isotropy 	$\Delta (S_{8}\times S_1^2)$.}
	\label{fig:gauss_a}
\end{figure}

Having established the invariance properties of the loss function, we now turn to investigate how critical points  of $\ploss_\star$ reflect this symmetry. The numerical results are obtained by first initializing the 
entries of $W$ to i.i.d. centered Gaussians with variance $1/d$, an initialization scheme often referred 
to as Xavier initialization. The vanilla gradient descent algorithm is then used to minimize (\ref{prob:opt}) until 
the gradient norm has been driven below a threshold of $1\mathrm{e-}6$ (unless otherwise stated). 
For convenience, we provide here the gradient expressions for the loss function in terms of the representation given in (\ref{loss_in_ker}) with $\ker_{\cN,3}$ replaced by a general kernel $\ker:\reals^d\times\reals^d\to\reals$ function. Assuming differentiability of the loss function at $W$,
\newcommand{\vh}{{h}}
\begin{align}
	\nabla \ploss_\ker(W) &= \sum_{i=1}^{k} \e_i\otimes \prn*{\sum_{j=1}^k  
		\kappa_\w
		(\w_i,\w_j)	- \sum_{j=1}^\vh \kappa_\w (\w_i,\v_j)},
\end{align}
where $\kappa_\w$ denotes the derivative of $\kappa$ with respect to the first argument, and $\e_i$ denotes the $i$-th unit vector. (The kernel representing optimization problem (\ref{prob:opt}) under the Frobenius norm is given by  $k(\cdot,\cdot)=\inner{\cdot,\cdot}^n$.) Upon convergence, the weight matrices are permuted so as to align with 
$\Delta (S_{d_1}\times S_{d_2}\times \cdots \times S_{d_p})$-isotropy group for $d_1\ge 
d_{2}\ge \dots\ge d_p$, when possible.  The Hessian expression, given for (\ref{loss_in_ker}), is,
\begin{align}
	\nabla^2 \ploss_\ker(W) 
	&= \sum_{i,j=1}^{k} \e_i \e_j^\top\otimes \kappa_{\w,\v} (\w_i,\w_j)	\\
	&+ \sum_{i=1}^{k} \e_i \e_i^\top\otimes 
	\prn*{\sum_{j=1}^k \kappa_{\w,\w} (\w_i,\w_j) - \sum_{j=1}^\vh \kappa_{\w,\w} 
		(\w_i,\v_j)},\nonumber
\end{align}
with $\kappa_{\w, \w}$ (resp. $\kappa_{\w, \v}$) denoting the derivative of $\kappa_\w$ with respect to $\w$ (resp. $\v$), assuming twice differentiability.\\
\\
\paragraph{Identity target tensor}\label{sec:identity}
First, we consider (\ref{prob:opt}) with $T=\sum_{i=1}^h \e_i$, assuming $d=k=h$ and $n\in\{3,5\}$. 
Spurious minima obtained by repeating the training procedure described earlier exhibit, consistently over all experiments, large isotropy subgroups of $S_d\times S_d$. Several minima obtained under the Frobenius norm are presented in \pref{fig:r3_a} and \pref{fig:r5_a}. Minima obtained under the cubic-Gaussian norm are presented in \pref{fig:gauss_a}. All are seen to be symmetry breaking.\\
\\

\begin{figure}[h]
	\centering 
	\includegraphics[scale=0.35]{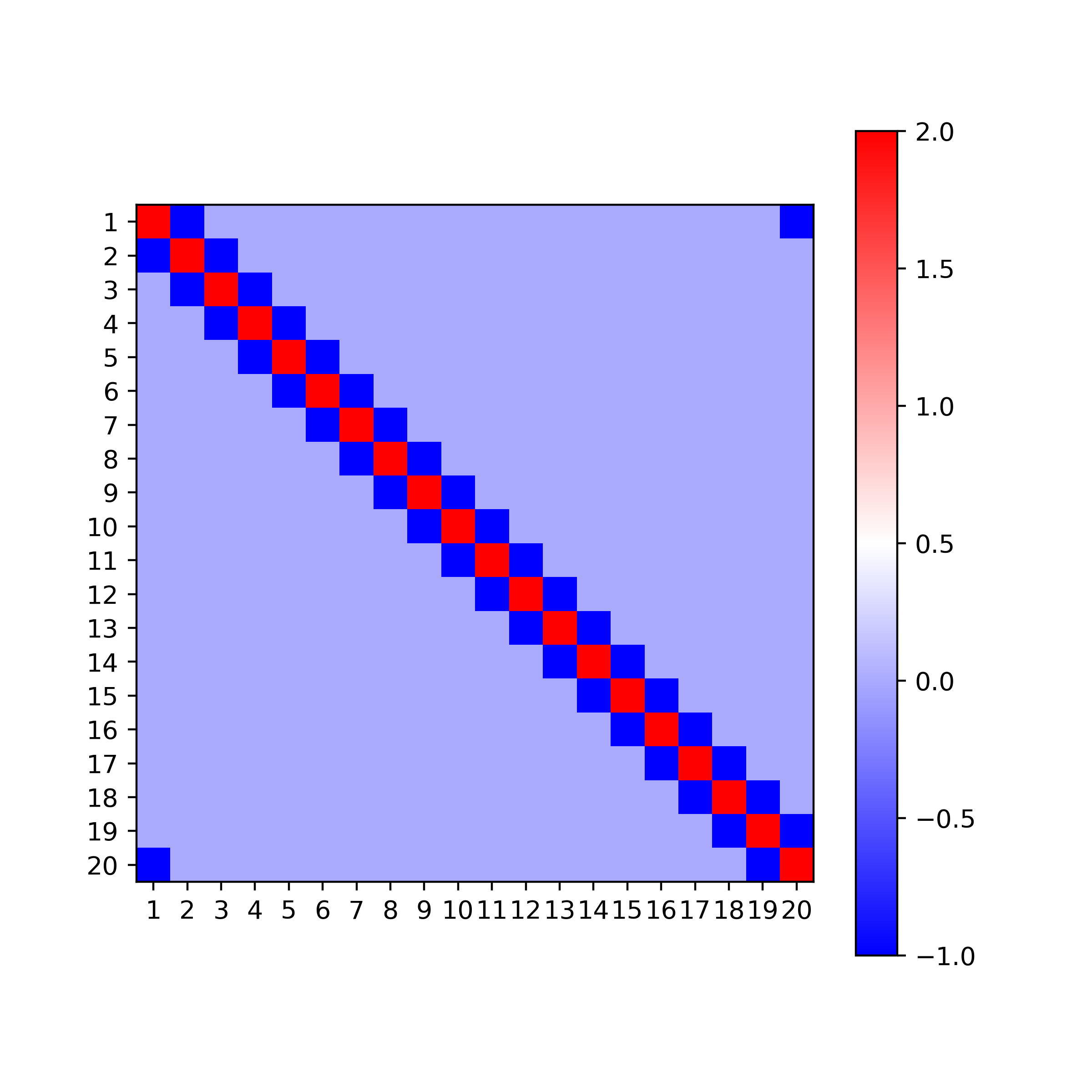}
	\includegraphics[scale=0.35]{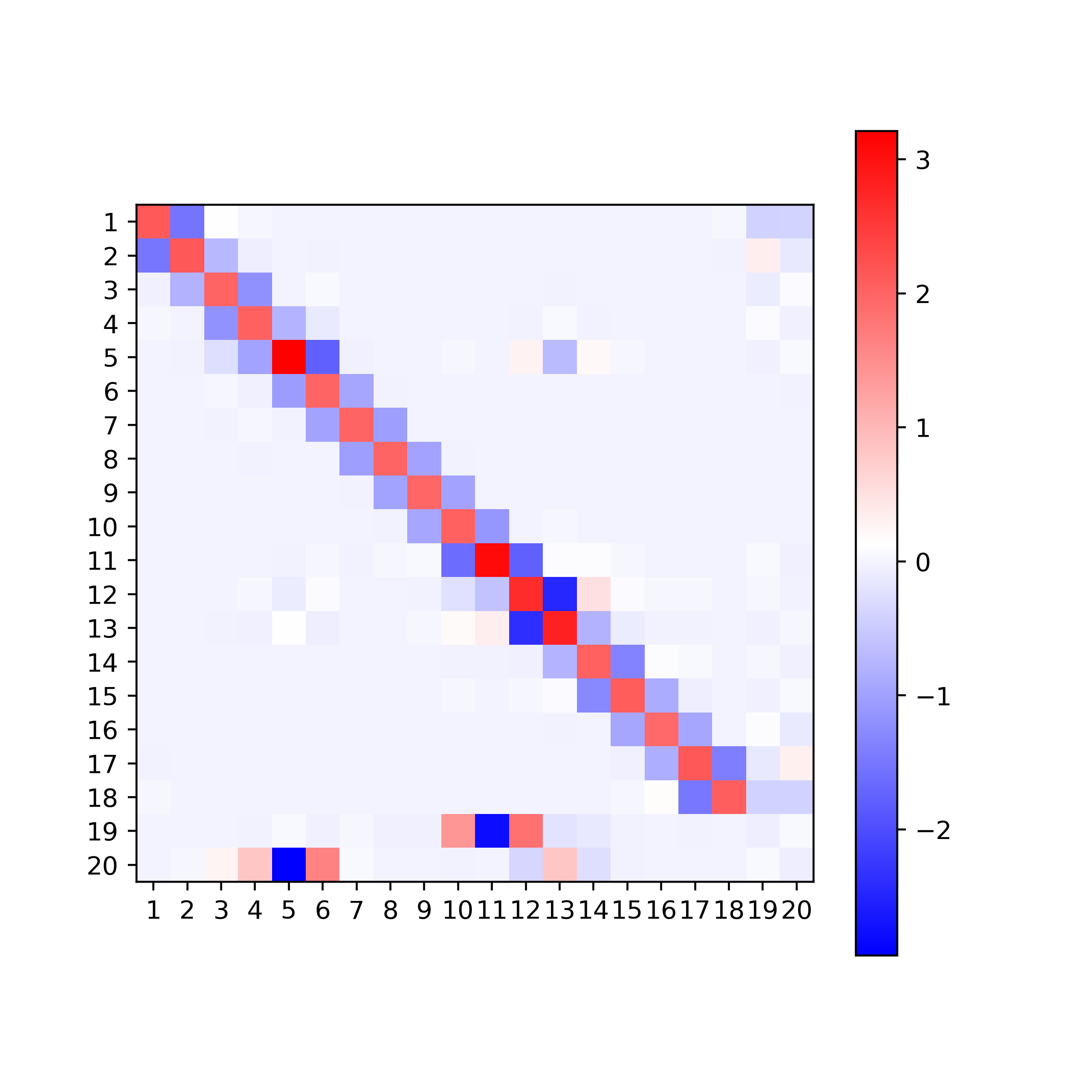}
	\caption{
		An approximated critical point of $\ploss_\cN$ for 
 $d=k=20$ and the Laplacian target tensor 
		$T_{\text{Lap}}$. (Left) the target 
		weight matrix $V_{\text{Lap}}$. (Right)	a critical point computed numerically (with gradient norm 		$<1\mathrm{e-}2$)}	\label{fig:conv}
\end{figure}

\paragraph{Convolutional target tensor}\label{sec:conv}
We now consider convolutional target tensors, demonstrating in particular symmetry breaking for non-identity target tensors. Specifically, we study problem (\ref{prob:opt}) with target tensor 
$T_{\text{Lap}} := \sum_{i=1}^d \v_i$,  where $\v_i$ denote the rows of the circulant \textit{Laplacian-filter} matrix
\begin{align}
	V_{\text{Lap}} := \begin{pmatrix}
		-2 & 1& 0 &\cdots & 0 & 1 \\\\
		&\ddots &&\ddots &&\\\\
		1 & 0 & \cdots & 0& 1 & -2 \\
	\end{pmatrix}.
\end{align}
The circulant structure of $V_{\text{Lap}}$ gives a $\Delta 
\bZ_d$-isotropy, corresponding to convolutional filters used in 
contemporary artificial neural network architectures. Spurious minima obtained by repeating the training procedure described earlier exhibit, consistently over all runs, large isotropy subgroups of $S_d\times \bZ_d$.  An example of an approximated critical point is given in \pref{fig:conv}.

\section{Conclusion}

In this note, we present a numerical study indicating that symmetry breaking phenomena occur for symmetric tensor decomposition problems. Our formulation of the invariance properties of the loss function (\ref{loss_in_ker}) subsumes a wide class of nonconvex kernel-like nonconvex optimization problems.
A particular implication is that symmetry breaking phenomena of the nature studied in this work applies more generally than for two layer ReLU networks, the setting in which the phenomenon was first observed. The numerical results  provided were later used in \cite{arjevani2023symmetry} to study and characterize infinite families of critical points.

\section*{Acknowledgments}
YA acknowledge partial support from the Israel Science Foundation (grant No. 724/22).
JB acknowledges partial support 
from the Alfred P. Sloan Foundation, NSF RI-1816753, NSF CAREER CIF 1845360, 
NSF CHS-1901091 and Samsung Electronics. JK acknowledges partial support 
from the Simons Collaboration in Algorithms and Geometry and start-up grants 
at UT Austin.

\bibliographystyle{ieeetr}
\bibliography{bib}

\end{document}